\documentclass[11pt,a4paper]{article} 

\usepackage{anysize}
\usepackage[format = hang]{caption}
\usepackage{amsmath,amsthm,verbatim,amssymb,amsfonts,amscd, graphicx}
\usepackage{graphics,natbib}
\usepackage{titlesec,footmisc}
\usepackage{hyperref} 

\usepackage{listings}
\usepackage{booktabs} 

\theoremstyle{plain}

\theoremstyle{definition}

\def\E{{\rm E}}
\def\Var{{\rm Var}}

\def\hatt{\widehat}
\def\arr{\rightarrow}

\def\half{\hbox{$1\over 2$}}
\def\eps{\varepsilon}

\def\beq{\begin{eqnarray}}
\def\eeq{\end{eqnarray}}

\def\beqn{\begin{eqnarray*}}  
\def\eeqn{\end{eqnarray*}}

\def\E{{\rm E}}
\def\Var{{\rm Var}}

\def\Pr{P}

\def\arr{\rightarrow}
\def\hatt{\widehat}

\def\sumin{\sum_{i=1}^n}

\def\eps{\varepsilon}
\def\half{\hbox{$1\over2$}}

\titleformat{\section}{\normalfont\large\sc\centering}{\thesection}{1em}{}
\titleformat{\subsection}[runin]{\normalfont\large\bfseries}{\thesubsection}{1em}{}
\numberwithin{equation}{section} 
\renewenvironment{abstract}
               {\list{}{\rightmargin\leftmargin}%
                \item[\text{\hspace{10mm}\sc Abstract.}]\relax}
               {\endlist}

\begin{document}

\def\heute{June 2000}

\begingroup
\begin{centering} 

  \Large{\bf A Note on Kernel Density Estimators \\ With Optimal Bandwidths}
  \\[0.8em]
\large{\bf Nils Lid Hjort$^{\rm a}$ and Stephen G.~Walker$^{\rm b}$} \\[0.3em] 

\smallskip 
\small {\sc $^{\rm a}$Department of Mathematics, University of Oslo, Norway} \\[0.3em]
\small {\sc $^{\rm b}$Department of Mathematical Sciences, University of Bath, UK}

\smallskip 
\small {\sc {\heute}}\par
\end{centering}
\endgroup


\begin{abstract}
\small{We show that the cumulative distribution function
corresponding to a kernel density estimator with optimal bandwidth
lies outside any confidence interval, around the empirical distribution
function, with probability tending to 1 as the sample size increases.

}
\end{abstract}


\section{Introduction} 

Kernel density estimation is a popular method for estimating the 
probability density
function (pdf) of an observed data set which obviates the need for a
parametric model. Much has been written on the subject as a consequence
of the problem of selecting a bandwidth. A recent review is provided
by Wand and Jones (1995).

In this paper we highlight a property of the 
classically recommended kernel density estimators (bandwidths of size
$O(n^{-1/5})$, where $n$ is the sample size).
The property is that the cumulative distribution function (cdf)
corresponding to the density estimator 
falls outside every reasonable confidence interval
or band of the empirical cdf, for every point, with probability tending to 1
as the sample size increases. 
A bandwidth of size $O(n^{-1/4})$ corrects this.   

It could be argued that the behaviour of the cdf is of minor
importance if interest is in estimating a pdf. However, a glance
at Parzen (1962), one of the pioneering papers in the field of density
estimation,
will convince the reader that kernel density
estimation was motivated by attempts to obtain gradients
from the empirical cdf. If the empirical cdf were differentiable
the pdf estimator would, according to Parzen at least, be the
pdf corresponding to this empirical cdf. Consequently, the connection between 
a kernel density estimator and the empirical cdf via the
cdf corresponding to the density estimator should not be
ignored. Even if one's thinking is solely on the density estimator
and optimal bandwidths, knowing that as the sample size tends
to infinity  the cdf of the density estimator leaves all confidence
intervals around the empirical cdf with probability going to 1, should
be a cause for concern. The result suggests the optimal bandwidth
is oversmoothing.

Let $X_1,\ldots,X_n$ be independent with common density function $f$ 
and cumulative distribution function $F$. 
The data give rise to their empirical cdf 
$F_n(x)=n^{-1}\sumin I\{X_i\le x\}$.
The classic nonparametric simultaneous confidence band
takes the form 
$${\rm CI}^{(1)}_n=[F_n(x)-c_n/\sqrt{n},
  F_n(x)+c_n/\sqrt{n}], \eqno(1.1)$$
where $c_n$ is the appropriate quantile of the distribution 
of $\sqrt{n}D_n$ where
\[D_n=\max_x|F_n(x)-F(x)|.\] 
In the limit as $n$ grows, $c_n$ becomes the quantile of the 
Kolmogorov--Smirnov distribution; for example, $c_n=1.224$ 
for 90\% confidence and $c_n=1.358$ for 95\% confidence. 
Alternatively, a confidence interval is based on the 
normal approximation to the binomial distribution,
$${\rm CI}^{(2)}_n=F_n(x)
  \pm d_n\{F_n(x)(1-F_n(x))\}^{1/2}/\sqrt{n}. \eqno(1.2)$$
With $d_n$ equal to 1.645 and 1.96 we have pointwise bands 
with approximate confidence level 90\% and 95\% for each $x$, 
while choosing $d_n$ equal to 2.89 and 3.15 provides global bands
with simultaneous confidence approximately 90\% and 95\%, 
valid for all $x$ between the 0.05 and 0.95 quantiles 
of the underlying $F(x)$ distribution.
The specific confidence interval is not relevant to our result.


For a symmetric kernel density function $k(u)$,
with associated cumulative distribution function $K(u)$, consider 
$$\hatt F_h(x)=n^{-1}\sumin K(h^{-1}(x-X_i)). $$
This is the smoothed empirical distribution function, 
inextricably linked to its more famous derivative, 
the kernel density estimator 
$$\hatt f_h(x)=n^{-1}\sumin h^{-1}k(h^{-1}(x-X_i)).$$ 
In Section 2 we state and prove the result.

\section{The result} 


Among the literature one of the strongest messages
is that $h$ should tend to zero with speed $n^{-1/5}$. See, for example,
Silverman (1986), Wand and Jones (1995). In this case, we prove the
following:

\smallskip
{\it Theorem.} With probability tending to 1 as $n$ tends to infinity,
the cdf of the  optimal kernel density estimator
will land outside all confidence bands around the empirical cdf,
including simultaneous and pointwise ones.

\smallskip
Proof. Consider the variable 
$$Z_n(x)=n^{1/2}\{\hatt F_h(x)-F_n(x)\}
  =n^{-1/2}\sumin A_i(x), \eqno(2.1)$$
writing $A_i(x)=K(h^{-1}(x-X_i))-I\{X_i\le x\}$. 
Saying that $Z_n(x)$ is outside $[-c_n,c_n]$ is
the same as stating that $\hatt F_h(x)$ lies outside the classic band (1.1).
The following calculations aim to find out what 
happens to $Z_n(x)$ as $n$ increases. 

First consider the mean. With substitution and partial integration, 
the mean of the first term of $A_i(x)$ can be written 
$$\int K(h^{-1}(x-y))f(y)\,\d y=\int K(v)hf(x-vh)\,\d v
  =\int k(v)F(x-hv)\,\d v. $$
A Taylor expansion gives $F(x)+\half k_2h^2f'(x)+o(h^2)$,
where $k_2=\int u^2k(u)\,\d u$. Hence 
$$\E Z_n(x)=\half k_2n^{1/2}h^2f'(x)+o(n^{1/2}h^2), $$
where the remainder term typically would be of size $O(n^{1/2}h^4)$
(requiring three derivatives of $f$ to exist at $x$). 
For $h=an^{-1/5}$, the recommended choice, we have 
$$\E Z_n(x)=\half k_2a^2n^{1/10}f'(x)+O(n^{-3/10}). \eqno(2.2)$$ 

Next, we square $A_i(x)$ and work with each of the terms separately. 
The first term required is
\beqn
\int K(h^{-1}(x-y))^2f(y)\,\d y
&=&\int K(v)^2hf(x-hv)\,\d v \\ 
&=&2\int K(v)k(v)F(x-hv)\,\d v
\eeqn 
(using here the fact that $K(v)^2F(x-hv)$ tends to zero at both ends).
The second necessary calculation is 
\beqn 
\int K(h^{-1}(x-y))I\{y\le x\}f(y)\,\d y
&=&\int_0^\infty K(v)hf(x-hv)\,\d v,  \\ 
&=&\half F(x)+\int_0^\infty k(v)F(x-hv)\,\d v.
\eeqn 
This leads to 
$$\E A_i(x)^2
=2\int K(v)k(v)F(x-hv)\,\d v-F(x)
  -2\int_0^\infty k(v)F(x-hv)\,\d v+F(x). $$
The leading terms of the variance of $Z_n(x)$ are accordingly
$$\Var\,A_i(x)=2(e_1-d_1)hf(x)-(e_2-d_2)h^2f'(x)+\hbox{$1\over 3$}(e_3-d_3)h^3f''(x)+O(h^4), \eqno(2.3)$$
where $e_j=\int_0^\infty v^jk(v)\,\d v$ 
and $d_j=\int_{-\infty}^\infty v^jk(v)K(v)\,\d v$. 
Now we can write $K(v)=1/2+vk(w)$ where $0<w<v$. Therefore,
$d_1=2\int_0^\infty v^2k(u)k(w)\,\d v$ which is less than $e_1$ since
$2vk(w)<1$. Using the same expansion of $K(v)$ we can prove that $e_2=d_2$.

So, with bandwidths chosen as $h=an^{-1/5}$, we have $EZ_n(x)=\alpha_n$,
with leading term $\half k_2a^2f'(x)n^{1/10}$, and
$\Var\,Z_n(x)=\beta_n$, with leading term $2(e_1-d_1)af(x)n^{-1/5}$.
Then $\pi_n=\Pr\{-c<Z_n<c\}=\Pr\{\alpha_n-c<T_n<\alpha_n+c\}$, 
where $T_n=\alpha_n-Z_n$, so clearly $\pi_n\rightarrow 0$, whenever 
$f'(x)\ne 0$.
 
Rephrasing, the heart of the matter is that $Z_n(x)$ climbs 
slowly towards plus or minus infinity, depending on the sign
and size of $f'(x)$. The calculation given shows that for each
single $x$ at which $f'(x)$ is nonzero, $\hatt F_h(x)$ 
eventually lands outside all natural confidence bands of
the types (1.1), (1.2) and so on. 
 
We ought to point out that the convergence towards 1 of 
not belonging to the confidence bands is quite slow. 

\section{Some Remedies and Further Results} 


\subsection{An asymptotic representation.}

We start this section by establishing limiting normality 
of $Z_n(x)$, properly normalised. This will lead to a 
useful asymptotic representation of $Z_n(x)$, and 
also make it possible to accurately compute the probability 
that the kernel cumulative estimator actually falls outside
the natural bands (as opposed to only demonstrating that 
the limiting probability is one). 

Consider
$$M_n(x)=\frac{Z_n(x)-\E Z_n(x)}{V^{1/2}f(x)^{1/2}h^{1/2}}
=\frac{\sumin [A_i(x)-\E A_i(x)]}{V^{1/2}f(x)^{1/2}h^{1/2}n^{1/2}},$$
where $V=2(e_1-d_1)$.
It has mean zero and variance 
of the form $\Var\,A_i(x)$ divided by $Vf(x)h$. 
Via equation (2.3) it is clear that the variance tends to one
if only $h\arr0$. Hence, by the Lindeberg and Feller theorems, 
$M_n(x)$ is asymptotically a standard normal if and only if 
$$\sumin \E\Big|{A_i(x)-\E A_i(x)\over V^{1/2}f(x)^{1/2}n^{1/2}h^{1/2}}\Big|^2
  I\Bigl\{\Big|{A_i(x)-\E A_i(x)\over V^{1/2}f(x)^{1/2}n^{1/2}h^{1/2}}\Big|
   \ge \eps\Bigr\}\arr0 $$
for each $\eps$. But since the $A_i(x)$ variables are i.i.d., 
this reduces to the requirement 
$$\E{(A_i(x)-\E A_i(x))^2\over Vf(x)h}
  I\bigl\{|A_i(x)-\E A_i(x)|\ge\eps V^{1/2}f(x)^{1/2}(nh)^{1/2}\bigr\}
  \arr 0. $$
And since $A_i(x)$ and its mean are bounded, by 1, this
Lindeberg condition holds whenever $nh\arr\infty$. 

Looking back to the earlier approximation of the mean of $Z_n(x)$,
let us also represent the process in the form 
$$Z_n(x)=\half k_2n^{1/2}h^2f'(x)+V^{1/2}f(x)^{1/2}h^{1/2}N_n(x). \eqno(3.1)$$
Here 
$$N_n(x)=M_n(x)+r_n(x)/\{V^{1/2}f(x)^{1/2}h^{1/2}\},$$ 
where $\E Z_n(x)=\half k_2n^{1/2}h^2f'(x)+r_n(x)$. 
Hence $N_n(x)$ is a limiting standard normal provided $r_n(x)=o(h^{1/2})$. 
When $f$ is smooth at $x$, $r_n(x)=O(n^{1/2}h^4)$,
so $N_n(x)$ in representation (3.1) is a standard normal 
in the limit provided $h\arr0$, $nh\arr\infty$ and $nh^7\arr0$.

Under these conditions, 
an approximation to the probability of belonging to the right set, 
as determined by (1.1), is 
$$\pi_n \approx\Pr\Bigl\{{-c-\half k_2f'(x)n^{1/2}h^2\over V^{1/2}f(x)^{1/2}h^{1/2}}
  \le N(0,1) \le {c-\half k_2f'(x)n^{1/2}h^2\over 
   V^{1/2}f(x)^{1/2}h^{1/2}}\Bigr\}.  \eqno(3.2)$$
If $f'(x)=0$ then $\pi_n\rightarrow 1$. If $f'(x)\ne 0$ then a
number of possibilities arise. 
For the interesting case of $h=an^{-1/4}$, we instead find that 
$$Z_n(x)=\half k_2a^2f'(x)+V^{1/2}a^{1/2}f(x)^{1/2}n^{-1/8}N_n, $$
which means that in the limit, inclusion in the  
interval (1.1) is 1 or 0, depending on whether $\half k_2a^2f'(x)$
is within $[-c,c]$ or not.
If it were known that $f'(x)>0$, then for inclusion we would
need
\[a<\sqrt{\frac{2c}{f'(x)k_2}}.\]
 In the case of $an^{-\eps}$, with $\eps<{1\over 4}$,
the limiting probability of $\hatt F_h(x)$ belonging to (1.1) is one.


\subsection{New bandwidth rules that do not oversmooth.} 

The above results suggest that the maximum rate with which 
$h$ should go to zero, or equivalently the maximum amount
of smoothing allowed by a data set of size $n$ in the purely 
nonparametric setting, is given by $h=an^{-1/4}$. 

A practical idea for a maximum smoothing parameter  
is based on the insistence that $\hatt F_h$ must belong 
to confidence band (1.1) for all $x$. This leads to suggesting 
$$\hatt h=\sup\{h>0\colon\hatt F_h(x)\in {\rm CI}^{(1)}_n
  {\rm\ for\ all\ }x\}=\sup\{h>0\colon \max_x|Z_n(x)|<c_n\} \eqno(3.3)$$
as the `maximum smoothing bandwidth'. Note that
$$\max_x|Z_n(x)|=\max_i\left\{\max\left\{|Z_n(X_i)|,|Z_n(X_i-)|\right\}
\right\},$$
so it is easy to compute $\hatt h$ for any given data set. A brief study
is carried out in Section 4.
One might suggest
selecting as final $h$ the one minimising the cross validation
curve subject to the constraint $h\le \hatt h$, for example.

To learn about the behaviour of $Z_n(x)$ and its maximum absolute
value, note from (2.1) that 
${\rm cov}\{Z_n(x),Z_n(y)\}={\rm cov}\{A_i(x),A_i(y)\}$, 
and the size of this covariance tends to zero as $h\arr0$
for each fixed pair $(x,y)$. This holds generally, but let 
us illustrate it for the case of $K$ supported on a bounded
interval, which we may take to be $[-\half,\half]$. 
Then $A_i(x)$ is always zero outside $x\pm\half h$,
so that the covariance in question is of size 
$-{1\over4}k_2^2h^4f'(x)f'(y)$ whenever $y$ is more than $h$ away from $x$;
in other words, the correlation between $Z_n(x)$ and $Z_n(y)$ 
becomes $O(f'(x)f'(y)h^3)$. Hence $M_n(x)$ and $N_n(x)$ in the 
above representations of $Z_n(x)$ behave for large $n$ as 
white noise with variance 1. 

It is clear from previous results that $\hatt h$ 
must tend to zero faster than $an^{-1/5}$. 
With $h=an^{-1/4}$ one sees from (3.1) that $\max|Z_n(x)|$
goes to $\half k_2a^2|f'(x_0)|$ when $n\arr\infty$,
where $x_0$ is a point at which $|f'|$ is maximal. 
A rough approximation to the parameter of maximal smoothing 
is therefore $an^{-1/4}$ where $a=(2c/k_2)^{1/2}\|f'\|^{-1/2}$
and $\|f'\|=\max_x|f'(x)|$. A normal-based quick rule 
would accordingly be $(2c/k_2)^{1/2}\phi(1)^{-1/2}\hatt\sigma/n^{1/4}$,
with $\hatt\sigma$ the standard deviation of the data,
and with $c$ from (1.1) dictated by the wished for 
confidence level in the Kolmogorov--Smirnov band. 

More careful approximations to $\max_x|Z_n(x)|$ may be put forward,
involving the $x_0'$ that maximises the exact absolute 
mean of $Z_n(x)$ rather than its approximation, and adding 
a factor times the standard deviation 
$V^{1/2}f(x_0')^{1/2}a^{1/2}n^{-1/8}$. 
It would however be besides our main point to overanalyse
the behaviour of the (3.3) quantity.

\section{Discussion} 


We have highlighted a surprising and perhaps embarrassing property
of the classically optimal kernel density estimators. 
Thinking nonparametrically, the benchmark is the empirical cdf
and the associated nonparametric confidence intervals or bands. 
These are all we have, and they should be respected. Density estimation 
is motivated by obtaining gradients from distribution functions.
Why are the optimal gradients coming from a cdf which, as the 
size of the data grows, lies with ever increasing probability outside
the confidence interval?
We believe the property of the cdf lying inside confidence intervals
should approach probability 1 as $n\rightarrow\infty$ and that this
property dominates all others. Consequently, the $O(n^{-1/5})$ rule,
which violates this, is not asymptotically optimal, regardless
of its other asymptotic properties. For this reason $O(n^{-1/4})$ rules
can be regarded as being asymptotically optimal. That is,
minimising asymptotic mean squared error subject to the cdf 
lying inside confidence bands with probability going to 1.
This is the practical relevance of the paper.  

At the heart of the conflict is that two different 
loss function are at work, respectively squared error 
for $f$ and squared error for $F$. As regards squared error
loss for $F$ the natural and hard-to-beat estimator is 
the empirical $F_n$ (it is the uniformly minimum-variance unbiased estimator),
while the $O(n^{-1/5})$ bandwidth is 
demonstrably optimal from the points of view of 
mean squared error and integrated mean squared error for $f$. 
Traditionally these results are phrased in a framework 
of asymptotics but one may argue that the $n^{-1/5}$
result is also valid for finite $n$; Glad, Hjort and Ushakov (1999) 
give an exact, tight upper bound for the mean integrated 
squared error, under minimal assmptions on the density, 
which again leads to a bandwidth of type $an^{-1/5}$. 

Another answer is that two slightly different sets of 
assumptions are employed, perhaps implicitly. 
The traditional analysis of density estimator behaviour 
assumes that $f$ has two smooth derivatives,
while using the $F_n$ estimator is the minimalistic 
nonparametric estimator, assuming nothing except 
exchangeability of the data. There could be different
responses to this; one argument would be that the 
kernel estimators, as fine-tuned by the theoretically
strongest methods, often turn out to be oversmooth,
seemingly reflecting more smoothness than the data 
can promise. This is an argument for smoothing less,
and the suggestion of (3.3) is one viable
method, well-motivated without any assumptions of 
density smoothness; it is truly nonparametric and automatic.  

The concern of smoothing too much has also been touched 
in the kernel method literature, see e.g.~Scott (1992, Section 6.5).
The `oversmoothing bandwidth' defined there is however
still of size $O(n^{-1/5})$, and barely larger than 
the normal-rule-of-thumb proposal $1.059\,\sigma/n^{1/5}$;
by the result above even this oversmoothing limit is too big
for comfort, when $n$ is very large.

In cases where the statistician really believes in smoothness
of the underlying phenomenon, she should perhaps exploit 
it also when making inference about $F$. In this light,
the traditional $F_n$ estimator and its accompanying bands 
(1.1)--(1.2) are too weak, and can be improved upon. 
A better estimator would be $\hatt F_h$, for suitable small $h$,
and confidence bands can be constructed via proper study
of the process 
$$U_n(x)=n^{1/2}\{\hatt F_{h_1}(x)-\half k_2{h_1}^2\hatt f'_{h_2}(x)
    -F(x)\}, $$
with one bandwidth $h_1$ to control smoothing in $\hatt F_h$
and another to give a good estimate of $f'(x)$. 
It may be shown that the $U_n$ is asymptotically 
a time-transformed Brownian motion, so that the band 
$$\hatt F_{h_1}(x)-\half k_2{h_1}^2\hatt f'_{h_2}(x)
  \pm c\{\hatt F_{h_1}(x)(1-\hatt F_{h_1}(x))\}^{1/2}/\sqrt{n} $$
contains the full underlying $F$ curve with the same
limiting probability as the classic band (1.1), with the same $c$. 

The reason why this actually is a little bit better than (1.1),
at least for very large $n$, is that 
\beqn
\E\{\hatt F_h(x)-F(x)\}^2
&=& n^{-1}F(x)(1-F(x)) \\
& &\qquad -2d_1hn^{-1}f(x)
+\hbox{$1\over4$}k_2^2h^4f'(x)^2+O(h^2n^{-1}+h^6),
\eeqn 
with $d_1$ a positive constant, given in Section 2. 
The best $h_1$ is of size $O(n^{-1/3})$ and 
the theoretically best $h_2$ of size $O(n^{-1/7})$. 

\section*{References} 

\begin{description}

\item

{\sc Glad, I.K., Hjort, N.L. \& Ushakov, N.G.} (1999). 
Upper bounds of mean integrated squared error for kernel 
density estimation.  Unpublished manuscript.

\item 
{\sc Parzen, E.} (1962). On the estimation of a probability density function
and the mode. {\it Annals of Mathematical Statistics} {\bf 33}, 1065-1076.

\item

{\sc Scott, D.W.} (1992). {\it Multivariate Density Estimation: Theory,
Practice and Vizualization}, Wiley, New York.

\item 
{\sc Silverman, B.W.} (1986). {\it Density Estimation for Statistics and
Data Analysis}, Chapman and Hall, London.

\item

{\sc Wand, M.P. \& Jones M.C.} (1995). {\it Kernel Smoothing},
Chapman and Hall, London. 

\end{description}
 

\end{document}